\documentclass[12pt,a4paper,twoside]{amsart} % -*- latex -*-

\usepackage{amsfonts}
\usepackage{amsmath}
\usepackage{amsthm}
\usepackage{amscd}
\usepackage{euscript}

\DeclareMathOperator{\Sec}{Sec}
\DeclareMathOperator{\Ima}{Im}
\DeclareMathOperator{\Jac}{Jac}
\DeclareMathOperator{\Aut}{Aut}
\DeclareMathOperator{\expdim}{expdim}

\newcommand{\LL}{{\EuScript L}}
\newcommand{\FF}{{\EuScript F}_{H}}

\newcommand{\OO}{{\EuScript O}}

\newcommand{\PP}{{\mathbb P}}
\newcommand{\GG}{{\mathbb G}}

\newcommand{\poh}{p^{(0)},\dots,p^{(h)}}

\newcommand{\Poh}{P^{(0)},\dots,P^{(h)}}
\newcommand{\uunon}{u_1,\dots,u_n}

\newcommand{\C}{\mathbb{C}}
\swapnumbers
\newtheorem{thm}{Theorem}[section]
\newtheorem{lemma}[thm]{Lemma}

\newtheorem{prop}[thm]{Proposition}
\theoremstyle{definition}
\newtheorem{defin}[thm]{Definition}

\theoremstyle{remark}
\newtheorem{remark}[thm]{Remark}

\begin{document}
\title{Grassmann defectivity \`a la Terracini} \author{Carla Dionisi
  and Claudio Fontanari}

\address{Carla Dionisi\\
 Dipartimento di Matematica Applicata
  ``G.~Sansone'' \\
  Via S.~Marta, 3\\
 I-50139, Firenze, Italy}
\email{dionisi@dma.unifi.it}

\address{Claudio Fontanari \\
  Scuola Normale Superiore \\
  Piazza dei Ca\-va\-lie\-ri, 7\\
  I-56126, Pisa, Italy}

\email{fontanar@cibs.sns.it}
\begin{abstract}
  
 This work is a modern revisitation of a classical paper by Alessandro
Terracini, going back to 1915, which suggests an elementary but powerful
method for  studing  Grassmann defective varieties. In particular, the
case of Veronese surfaces is completely understood, giving an positive
answer to the so-called Waring problem for pairs of homogeneous polynomials
in three variables.

\vspace{0.3cm}
\noindent{\it Mathematics Subject Classification (2000)}: 14N15, 51N35
\end{abstract}
\maketitle 

%%%%%%%%%%%%%%%%%%%%%%%%%%%%%%%%%%%%%%%%%%%%%%%%%%%%%%%%%
\section{Introduction}
Let $V \subset \PP^r$ be an irreducible nondegenerate projective
variety of dimension $n$ defined over the complex field $\C$.
\begin{defin} 
  The \emph{$h$-secant variety} $\Sec_h(V)$ of $V$ is the irreducible
  variety given by the Zariski closure of the set:

\noindent$\{ p \in\PP^r$ s.~t. $p$ lies in the span of $h+1$ independent
points of $V$\}.
\end{defin}
\noindent Just counting parameters we get 
\begin{equation}
\label{eq:sec_hdim}
\dim \Sec_h(V) \leq   \min \{(n+1)(h+1)-1, r\}  
\end{equation}
where the right hand side is called the {\it expected dimension of
  $\Sec_h(V)$}. If strict inequality holds in (\ref{eq:sec_hdim}),
then $V$ is said to be \emph{$h$-defective} and the positive integer
\[
\delta_h(V):= \min \{(n+1)(h+1)-1, r\} - \dim \Sec_h(V)
\]
is called the \emph{$h$-defect} of $V$.

A basic tool for understanding defective varieties is the classical
Terracini's lemma (see \cite{Terracini:11} for the original version
and \cite{Dale:84} for a modern proof), which says that, for
$P_1,\dots, P_{h+1}$ general points and $P$ general in their span one
has
\begin{equation*}
  T_P(Sec_h(V)) = < T_{P_1}(V),\dots, T_{P_{h+1}}(V)>.
\end{equation*}
Therefore $V$ is $h$-defective if and only if 
\begin{equation*}
  \dim (< T_{P_1}(V),\dots, T_{P_{h+1}}(V)>) <  \min \{(n+1)(h+1)-1, r\}.
\end{equation*}

The systematic study of defective varieties goes back to the old
Italian school: we wish to mention at least the contributions of
Francesco Palatini (\cite{Palatini:06}, \cite{Palatini:09}), Gaetano
Scorza (\cite{Scorza:08}, \cite{Scorza:08bis}, \cite{Scorza:09}) and
Alessandro Terracini (\cite {Terracini:11}, \cite{Terracini:21}). This
great amount of work was recently rediscovered by various authors,
among whom Luca Chiantini and Ciro Ciliberto; we refer to their papers
\cite{ChiCil:01} and \cite{ChiCil:01bis} for rigorous proofs and
powerful generalizations of the classical results in the field.

In the present paper, instead, we focus on another kind of
defectivity, the so-called Grassmann-defectivity.
\begin{defin} 
  The \emph{$(k,h)$-Grassmann secant variety} $\Sec_{k,h}(V)$ of $V$
  is the Zariski closure of the set:

\noindent$\{ l \in \GG(k,r)$ s.~t. $l$ lies in the span of $h\!+\!1$
independent points of $V\}.$
\end{defin}

As above, we have an obvious inequality
\begin{equation}
\label{eq:sec_k,hdim}
\dim \Sec_{k,h}(V) \leq \min \{ (h+1)n+(k+1)(h-k), (k+1)(r-k) \}   
\end{equation}
and we may introduce in a natural way the definitions of
\emph{expected dimension} of $ \Sec_{k,h}(V)$,
\emph{$(k,h)$-defectivity} and \emph{$(k,h)$-defect} $\delta_{k,h}(V)$
of $V$.

Unfortunately there seems to be no easy form of a Terracini type lemma
which may help in this situation. As a consequence, the problem of
classifying Grassmann-defective varieties is rather hard (see
\cite{ChiCop:01} for a first step in this direction; moreover, 
\cite{ChiCil:01ter} shows that there are no Grassmann defective curves). 
However, in the memoir
\cite{Terracini:15}, going back to 1915, Terracini suggests that the
condition of $(h,k)$-defectivity for a variety of dimension $n$ may be
translated into the condition of $h$-defectivity for a variety of
dimension $n+k$. More precisely, we have the following statement (to
be proved in section~\ref{proofs}):
\begin{prop}\label{segre}
Let $V \subset \PP^r$ be an irreducible nondegenerate projective
variety of dimension $n$.  Let $\sigma: \PP^k \times V \to \PP^{r(k+1)+k}$ 
be the Segre embedding of $\PP^k \times V$.
Then 
$V$ is $(k,h)$-defective with defect $\delta_{k,h}(V)=\delta$ if and only
if 
$\sigma(\PP^k \times V)$ is $h$-defective with defect
$\delta_h(\sigma(\PP^k 
\times V))=\delta$.
\end{prop}

This fact was pointed out by Terracini (see p.~97 of
\cite{Terracini:15}) only for Veronese surfaces, but it turns out to
hold in complete generality with almost the same proof. The interest
of Terracini in the case of Veronese varieties is explained by the
simple observation that the Veronese variety $V_{n,d}$ of dimension
$n$ and degree $d$ is \emph{not} $(k,h)$-defective if and only if the
following Waring type problem admits an affirmative answer:

\vspace{5pt}
\noindent {\it Given positive integers $d$, $n$, $k$, $h$, may we write
  any $(k+1)$ homogeneous polynomials $f_j(x_0,\dots, x_n)$,
  $j=0,\dots k$, of degree $d$ as linear combinations of the same
  $(h+1)$ $d$-th powers of linear forms $l_i(x_0,\dots,x_n)$,
  $i=0,\dots,h$ ? }

\vspace{5pt}

Along these lines, Terracini's approach leads to the following
conclusions:
\begin{thm}\label{veronese}
If $(d,h)\neq(3,4)$ then $V_{2,d}$ is \emph{not} $(1,h)$-defective.
\end{thm}

Indeed the paper \cite{Terracini:15} contains a proof of this result under
the 
additional numerical hypothesis:
\[
h+1 \geq \frac{(d+1)(d+2)}{4}
\]
which can be removed using a subtler argument. 

We stress moreover that Theorem~\ref{veronese} is sharp, since the
Veronese surface $V_{2,3}$ is $(1,4)$-defective, as already
noticed by London in \cite{London:90} (see Remark \ref{last_rmk}).\\

This paper collects the results obtained by the authors under the
direction of Ciro Ciliberto during the EAGER Summer School PRAGMATIC
2001.  We would like to warmly thank Professor Ciliberto for his
patient guidance and his generous suggestions and the organizers of
PRAGMATIC (especially Professor Ragusa) for our nice stay in Catania.

This research was partially supported by MURST (Italy).

%%%%%%%%%%%%%%%%%%%%%%%%%%%%%%%%%%%%%%%%%%%%%%%%%%%%%%%%%%%%%%
\section{The proofs}\label{proofs}
\noindent 
\textbf{Proof of Proposition~\ref{segre}.} Consider the natural map
\begin{equation*}
  \phi: \underbrace{ \PP^h \times \dots \times \PP^h}_{k+1} \times V^{h+1} \to (\PP^r)^{k+1}
\end{equation*}
\[
((\lambda_{ij})_{\begin{tiny}  
  {\begin{split}
    i=0,\dots,k \\
    j=0,\dots,h
\end{split}}
\end{tiny}} , (p^{(j)})_{j=0,\dots,h})\mapsto (\sum_{j=0}^h \lambda_{ij}
p^{(j)})_{i=0,\dots,k}
\]
which to a set of coefficients and a $(h+1)$-ple of points lying on
$V$ associates a $(k+1)$-ple of points in $\PP^r$ each one contained
in the linear span of the $h+1$ points we started with. 

\noindent By definition, $V$ is
\emph{not} $(k,h)$-defective if and only if $\dim \Ima(\phi)$ has
maximal dimension.  If $p: U \subseteq \C^n \to V$ is a local
parametrization of $V$, we may introduce $\Jac(\phi \circ p)$, the
Jacobian matrix of $\phi \circ p$, given by
\begin{eqnarray*}
\begin{pmatrix}
\frac{\partial (\phi \circ p)}{\partial \lambda_{00}}
\dots \frac{\partial (\phi \circ p)}{\partial \lambda_{k0}} 
& \dots \dots&
\frac{\partial (\phi \circ p)}{\partial \lambda_{0h}}
\dots \frac{\partial (\phi \circ p)}{\partial \lambda_{kh}}& 
  \frac{\partial (\phi \circ p)}{\partial u_{1}^{(0)}} \dots
  \frac{\partial (\phi \circ p)}{\partial u_{n}^{(0)}}& \dots&
  \frac{\partial (\phi \circ p)}{\partial u_{1}^{(h)}} \dots 
  \frac{\partial (\phi \circ p)}{\partial u_{n}^{(h)}}
\end{pmatrix} = \\
\begin{pmatrix}
 p^{(0)} \ 0 \  \dots \ 0 & \dots \dots & p^{(h)}\ 0 \ \dots \ 0 &
      \lambda_{00}p_{u_1}^{(0)} \dots \lambda_{00}p_{u_n}^{(0)}& \dots
&\lambda_{0h} p_{u_1}^{(h)} \dots \lambda_{0h}p_{u_n}^{(h)} \\
      0 \ \ p^{(0)} \ \dots \ 0 & \dots \dots & 0 \ \ p^{(h)}\ \dots \
 0 &  \lambda_{10}p_{u_1}^{(0)} \dots \lambda_{10}p_{u_n}^{(0)}& \dots &
\lambda_{1h}p_{u_1}^{(h)} \dots \lambda_{1h}p_{u_n}^{(h)}\\
      \vdots& \vdots&  \vdots& \vdots& \vdots& \vdots \\
     \  0 \ \ \ 0 \ \dots \  p^{(0)} & \dots \dots & \ 0 \ \ \ 0 \  \dots
\  p^{(h)}
 & \lambda_{k0}p_{u_1}^{(0)} \dots
      \lambda_{k0}p_{u_n}^{(0)}& \dots & \lambda_{kh}p_{u_1}^{(h)}
      \dots \lambda_{kh}p_{u_n}^{(h)} \\
 \end{pmatrix}
\end{eqnarray*}

\noindent
and apply the inverse function theorem to conclude that $V$ is \emph{not} 
$(k,h)$-defective (resp., $(h,k)$-defective with defect $\delta$) if and
only 
if $\Jac(\phi \circ p)$ has maximal rank (resp., rank equal to the maximal
one 
minus $\delta$) at a general point.

On the other hand, consider the Segre embedding 
\begin{eqnarray*}
\sigma: \PP^k \times V & \to & \PP^{r(k+1)+k} \\ 
(\lambda_0, \dots, \lambda_k, p) &\mapsto& (\lambda_0 p, \dots, \lambda_k
p).
\end{eqnarray*}
If $t_i=\frac{\lambda_i}{\lambda_0}$, locally we have 
\begin{equation*}
(\sigma \circ p) (t_1,\ldots, t_k, \uunon)=(p(\uunon) ,t_1p(\uunon),\dots,
t_kp(\uunon))
\end{equation*}

and we may compute:
\begin{equation*}
  \begin{aligned}
    \frac{\partial (\sigma \circ p)}{\partial t_i}&=(0,\dots,0,p,0
\dots,0) 
&\quad \text{for } i=1,\dots k \\
    \frac{\partial (\sigma \circ p)}{\partial u_j}&=(p_{u_j},
t_1p_{u_j},\dots, 
t_kp_{u_j})= 
(\lambda_0 p_{u_j},\dots, \lambda_k p_{u_j}) &\quad \text{for } 
j=1,\dots n. \\
\end{aligned}
\end{equation*}
Hence $T(\sigma(\PP^k \times V))$ is spanned by the columns of the matrix 
\begin{equation*}
\begin{aligned}  
\quad &
\begin{pmatrix}
    \lambda_0 p & 0 & \dots & 0&\  \lambda_0 p_{u_1}& \dots &\lambda_0
p_{u_n} \\
    \lambda_1 p & p & \dots & 0&\  \lambda_1 p_{u_1}& \dots &\lambda_1
p_{u_n} \\
    \vdots           & \vdots &\vdots &\vdots &\vdots &\vdots &\vdots \\
    \lambda_k p & 0 & \dots & p&\  \lambda_k p_{u_1}& \dots &\lambda_k
p_{u_n} \\
\end{pmatrix}\\
\sim& \  
   \begin{pmatrix}
    \ \  p & 0 & \dots & 0&\  \lambda_0 p_{u_1}& \dots &\lambda_0 p_{u_n}
\\
    \ \  0 & p &\dots & 0&\  \lambda_1 p_{u_1}& \dots &\lambda_1 p_{u_n}
\\
    \vdots           & \vdots &\vdots &\vdots &\vdots &\vdots &\vdots \\
    \ \  0& 0 & \dots & p&\  \lambda_k p_{u_1}& \dots &\lambda_k p_{u_n}
\\
\end{pmatrix}
\end{aligned}
\end{equation*}
where the relation $\sim$ between these two matrices means that
their columns span the same vector space $A(\lambda_0, \dots, \lambda_k, p)$.

By Terracini's lemma,
\[
  T(\Sec_h(\sigma(\PP^k \times V))=<
A(\lambda_{00},\dots,\lambda_{k0},p^{(0)})
\vert \dots \vert  A(\lambda_{0h},\dots,\lambda_{kh},p^{(h)})>
\]
and $\sigma(\PP^k \times V)$ is \emph{not} $h$-defective (resp.,
$h$-defective with 
defect $\delta$) if and only if this matrix has maximal rank (resp., rank
equal to 
the maximal one minus $\delta$).  

Since $\Jac(\phi \circ p) \sim <
A(\lambda_{00},\dots,\lambda_{k0},p^{(0)}) \vert \dots \vert
A(\lambda_{0h}, \dots, \lambda_{kh},p^{(h)})>$ the thesis follows.
\qed

\begin{lemma} \label{hyperplane}
  Let $V \subset \PP^r$ be an irreducible nondegenerate projective
  variety of dimension $n$.  Let $\sigma: \PP^k \times V \to
  \PP^{r(k+1)+k}$ be the Segre embedding of $\PP^k \times V$.  Fix
  $\poh$ general points on $V$ and $\lambda^{(0)}, \ldots,
  \lambda^{(h)}$ general points in $\PP^k$, so that
  $P^{(j)}:=(\lambda_0^{(j)} p^{(j)},\ldots, \lambda_k^{(j)} p^{(j)})$
  is a general point on $\sigma(\PP^k \times V) \subset
  \PP^{r(k+1)+k}$ for $j= 0 \ldots h$; finally, take a general point
  $P \in < \Poh >$.  Then there is a natural identification between:
\begin{itemize}
\item hyperplanes $H \subset \PP^{r(k+1)+k}$ such that
  $T_P(\Sec_h(\sigma (\PP^k \times V))) \subset H$;
\item $k$-dimensional linear systems $\mathcal{H}$ of hyperplane
  sections of $V \subset \PP^r$ with a projectivity $\omega:
  \mathcal{H} \to \PP^k$ such that all the elements of the linear
  system pass through the points $p^{(j)} \in V$ and for every $j$ the
  hyperplane section of the linear system corresponding to
  $\lambda^{(j)}$ is tangent to $V$ at $p^{(j)}$.
\end{itemize} 
\end{lemma}

\begin{proof}
  If $X_{\alpha \beta}$ ($\alpha=0 \ldots k$, $\beta=0 \ldots r$) are
  coordinates on $\PP^{r(k+1)+k}$, the points of $\sigma(\PP^k \times
  V)$ are exactly those of the form
  \[X_{\alpha \beta} = \lambda_{\alpha}p(\uunon)_{\beta}\]
  in the
  notation of the proof of Proposition~\ref{segre}.  A hyperplane $H$
  in $\PP^{r(k+1)+k}$ is given by a linear equation
  \[\sum_{\alpha,\beta} a_{\alpha \beta} X_{\alpha \beta} = 0\]
  and if
  $H$ is general we may assume that $a_{\alpha \beta} \ne 0$ for every
  $\alpha, \beta$. So a hyperplane section of the variety
  $\sigma(\PP^k \times V)$ is of the form
\begin{equation*}
  \sum_{\alpha,\beta}a_{\alpha \beta} \lambda_{\alpha}p(\uunon)_{\beta}=0
\end{equation*}
i.e.
\begin{equation*}
 \sum_{\alpha} \lambda_{\alpha} \sum_{\beta} a_{\alpha
\beta}p(\uunon)_{\beta}=0.   
\end{equation*}
Hence every general hyperplane section of the variety $\sigma(\PP^k
\times V)$ corresponds to a $k$-dimensional linear system $\mathcal{H}$ 
of hyperplane sections of the variety $V$ with a fixed projectivity 
$\omega: \mathcal{H} \to \PP^k$; conversely, the data of such a
linear system $\mathcal{H}$ and a projectivity $\omega: \mathcal{H}
\to \PP^k$ uniquely determine a hyperplane section of $\sigma(\PP^k
\times V)$.  Moreover, if
\begin{equation*}
T_P(\Sec_h(\sigma(\PP^k \times V))) = < T_{P^{(0)}}(\sigma(\PP^k \times
V)), 
\ldots, T_{P^{(h)}}(\sigma(\PP^k \times V)) > \subset H
\end{equation*} 
then 
\begin{equation*}
  \begin{aligned}
    \sum_{\beta}& a_{\alpha \beta}p^{(j)}(\uunon)_{\beta}=0  &\quad
\forall \alpha=0,\dots, k \\
    \sum_{\alpha}& \lambda_\alpha^{(j)} \sum_{\beta}a_{\alpha
      \beta}{p_{u_{\gamma}}}^{(j)}_{\beta}=0 &\quad \forall
    \gamma=1,\dots, n.
\end{aligned}
\end{equation*}
In other words, all the elements of the linear system pass through the
points $p^{(j)} \in V$ and for every $j$ the hyperplane section of the
linear system corresponding to the values $\{ \lambda_\alpha^{(j)}
\}_{\alpha = 0 \ldots k}$ of parameters is tangent to $V$ at
$p^{(j)}$.
\end{proof}

\noindent
\textbf{Proof of Theorem~\ref{veronese}.} Assume by contradiction that
$V_{2,d}$ is $(1,h)$-defective. By Proposition~\ref{segre},
$\sigma(\PP^1 \times V_{2,d}) \subset \PP^{d(d+3)+1}$ has to be
$h$-defective, i. e.
\begin{equation}
\label{eqn:estimate}
\dim \Sec_h(\sigma(\PP^1 \times V_{2,d})) = d(d+3)+1- \delta < 4(h+1)-1
\end{equation}
with $\delta \ge 1$.

\noindent Hence if we take a general point $P$ as in the statement 
of Lemma~\ref{hyperplane}, then $T_P(\Sec_h(\sigma(\PP^1 \times
V_{2,d})))$ is contained in $\delta$ independent hyperplanes $H_t$ ($1
\le t \le \delta$), which we may assume to be general since $P$ is
general.

\noindent By Lemma~\ref{hyperplane}, each 
$H_t$ gives rise to a pencil ${\FF}_t$ of plane curves of degree $d$
all passing through $h+1$ general points $\poh$ in such a way that for
every $p^{(j)}$ at least one curve of the pencil passes doubly through
$p^{(j)}$.

\noindent By an infinitesimal Lemma already known to Terracini and 
reproved in modern times by Ciliberto and Hirschowitz in
\cite{CilHir:91}, we have that every $H_t$ is tangent to $\sigma(\PP^1
\times V_{2,d})$ along a positive dimensional variety $C_t$ passing
through $\Poh$. The points of $C_t$ are indeed Segre images of pairs
$(\lambda, p)$ such that the plane curve of ${\FF}_t$ corresponding to
$\lambda$ has a singular point in $p$. If  all the
$C_t$ are one-dimensional for each $H_t$ we have \emph{a priori} two cases:

(i) $\sigma^{-1}(C_t) = \bigcup_{j=0}^{h} \PP^1 \times \{p^{(j)} \}$,
so that all the curves of ${\FF}_t$ pass doubly through $\poh$;

(ii) $\sigma^{-1}(C_t)$ surjects on $\PP^1$ and injectively projects
to $\PP^2$ over a plane curve passing through $\poh$, so that all the
curves of ${\FF}_t$ have a double point.
 
\noindent Since the hyperplanes $H_t$ are general, by symmetry  the
same case occurs for all of them. Moreover, if the $C_t$'s are higher dimensional we fall
\emph{a fortiori} in case (ii).

In case (i) we have
\[{\FF}_t \subseteq \LL_{2,d}(2^{h+1}) \]
for every $t$, where
$\LL_{2,d}(2^{h+1})$ as usual denotes the linear system of plane
curves of degree $d$ with $h+1$ assigned general double points.  We
claim that
\[\dim <{\FF}_1 \ldots {\FF}_\delta > \ \ge \quad \frac{\delta}{2}.\]
To check the claim, let $ < {\FF}_1 \ldots
{\FF}_\delta >$ be spanned by the columns of a matrix
\begin{equation*}
  \begin{pmatrix}
    f_{01}& \ldots & f_{0 \delta}\\
    f_{11}& \ldots & f_{1 \delta}
  \end{pmatrix}
\end{equation*}
Just making elementary operations on columns, we obtain (up to
reindexing)
\begin{equation*}
\begin{pmatrix}
  f_{01}& \ldots & f_{0 \delta}\\
  f_{11}& \ldots & f_{1 \delta}
\end{pmatrix}
\sim 
\begin{pmatrix}
  f_{01}& \ldots  & f_{0x} & 0 & \ldots & 0\\
  f_{11}& \ldots  & f_{1x} & g_{1 x+1} &\ldots & g_{1 \delta}\\
\end{pmatrix}
\end{equation*}
where
\begin{equation*}
x = \dim <f_{01}  \ldots f_{0 \delta}>.
\end{equation*}

Since both $f_{01} \ldots f_{0x}$ and $g_{1x+1} \ldots g_{1 \delta}$
are linearly independent, we deduce
\[\dim < {\FF}_1 \ldots {\FF}_\delta > \  \ge \max \{x, \delta -x
\} \ge \frac{\delta}{2}\]
and the claim is checked.  \newline
\noindent By the claim, to get a contradiction it will be sufficient
to show that $ \dim \LL_{2,d}(2^{h+1}) < \frac{\delta}{2}$. Using
(\ref{eqn:estimate}) we compute:
\begin{equation*}
  \begin{split}
    \expdim \LL_{2,d}(2^{h+1}) &= \frac{d(d+3)}{2}-3(h+1) \le
    \frac{d(d+3)}{2}
    - \frac{3}{4}(d(d+3)+3-\delta)=  \\
    &= \frac{-d(d+3)-9+3 \delta}{4} < \frac{\delta}{2}.
 \end{split}
\end{equation*}

In order to conclude, just notice that if $d \ge 5$ then
$\LL_{2,d}(2^{h+1})$ is nonspecial by the Alexander-Hirschowitz
theorem (see \cite{AleHir:95}), so
\[\dim  \LL_{2,d}(2^{h+1}) = \expdim  \LL_{2,d}(2^{h+1}) <
\frac{\delta}{2};\]
if instead $d\leq4$, the only special system arises when $d=2$ and
$h=1$ or $d=3$ and $h=4$, but also in 
these cases we have (see \cite{Cil:01}, Example 4.3)
\[\dim  \LL_{2,d}(2^{h+1}) = 0 < \frac{\delta}{2}.\]
So case (i) is over. 

In case (ii), it follows from Bertini's theorem that all the pencils
have a base curve of the same degree, say $m$.  Since this curve has
to pass through $\poh$, we have
\begin{equation}\label{eqn:first}
\frac{m(m+3)}{2} \ge h+1  \ge \frac{d(d+3)+3-\delta}{4}.
\end{equation}
If the base curve pass doubly through $\poh$, we may argue exactly as in case 
(i); otherwise, the moving parts of the $\delta$ independent pencils
have to pass through $\poh$ in correspondence with prescribed general 
coefficients; hence the generic fiber of the natural map
\[\GG(\PP^1, \PP H^0(\PP^2, \OO_{\PP^2}(d-m))) \longrightarrow
(\PP^1)^{h+1} / \Aut(\PP^1),\]
which associates to a one-dimensional
linear system the $(h+1)$-ple of coefficients corresponding to its
curves through $\poh$, must have dimension $\ge \delta -1$. It follows
that
\[
\dim \GG(\PP^1, \PP H^0(\PP^2, \OO_{\PP^2}(d-m))) - (h+1) + \#\Aut(\PP^1) 
\ge \delta - 1, 
\]
i.e.
\begin{equation}\label{eqn:second}
(d-m)(d-m+3)+1 \ge \frac{d(d+3)+3- \delta}{4} + \delta -1.
\end{equation} 
Summing equations (\ref{eqn:first}) and (\ref{eqn:second}) we get 
\[
3 \frac{m(m+3)}{2}+(d-m)(d-m+3)+1 \ge d(d+3)+2
\] 
and we may deduce that 
\[
d-m < \frac{d+3}{5}.
\]
Substituting in (\ref{eqn:second}) we obtain 
\[
\frac{(d+3)(d+18)}{25}+1 > \frac{d(d+3)+2}{4}
\] 
hence  $d \leq 3$.

\noindent If $d=2$ then $m=1$, so (\ref{eqn:first}) implies $\delta
\geq 5$, while (\ref{eqn:second}) forces $\delta \leq \frac{11}{3}$,
contradiction.

\noindent If $d=3$ and $m=1$, then (\ref{eqn:first}) implies $\delta
\geq 13$, while  (\ref{eqn:second}) forces $\delta \leq 9$, contradiction.  
 
\noindent If $d=3$ and $m=2$, then (\ref{eqn:first}) implies $h \leq
4$. Hence for $(d,h)\neq(3,4)$ (\ref{eqn:first}) implies $\delta \geq 5$
 while  (\ref{eqn:second}) forces $\delta \leq 1$, contradiction.  

\noindent So the proof is over.
\qed

\begin{remark}\label{last_rmk}
  To check that $V_{2,3}$ is indeed $(1,4)$-defective we may argue as
  in the proof of Theorem~\ref{veronese}. In fact, since for $m=2$,
  $h=4$, $d = 3$ and $\delta = 1$ conditions (\ref{eqn:first}) and 
  (\ref{eqn:second}) are verified,
  there exists a pencil of plane cubics with a fixed component of
  degree $2$ passing through $5$ general points and a moving part
  passing through the same points in correspondence with prescribed
  general coefficients. Hence by Lemma~\ref{hyperplane} we deduce that
  $T_P(\Sec_4(\sigma(\PP^1 \times V_{2,3})))$ is contained in at least
  one hyperplane of $\PP^{19}$, so $\dim \Sec_4(\sigma (\PP^1 \times
  V_{2,3})) \le 18$. Since $\expdim \Sec_4(\sigma(\PP^1 \times
  V_{2,3}))=19$, it turns out that $\sigma(\PP^1 \times S)$ is
  $4$-defective. Now the thesis directly follows from
  Proposition~\ref{segre}.
\end{remark}
%%%%%%%%%%%%%%%%%%%%%%%%%%%%%%%%%%%%%%%%%%%%%%%%%%%%%%%%%

\end{document}